\documentclass[12pt]{article}
\input{amssym.def}
\input{amssym}

\makeatletter
\@addtoreset{equation}{section}
\makeatother

\def\be{\begin{equation}}
\def\ee{\end{equation}}
\def\ba{\begin{eqnarray}}
\def\ea{\end{eqnarray}}
\def\lb{\label}

\def\Tr#1{{\rm Tr\str{-1.3}}_{R^{\mbox{\scriptsize
$(#1)$}}}}
\def\str#1{\rule[#1mm]{0pt}{1mm}}

\def\tbl#1#2{{\ifmmode
\Bigl\{\!\!\begin{array}{c}
\scriptstyle #1\\[-2pt]
\raisebox{2pt}{$\scriptstyle #2$} \end{array}\!\!\Bigr\}
\else
\raisebox{2pt}{\scriptsize$\left\{\!\!\!\begin{array}{c}
#1\\[-1pt] #2 \end{array}\!\!\!\right\}$}\fi}}

\def\mtbl#1#2#3{{\ell_{#3}\hspace{-3pt}
\left\{\hspace{-1pt}
\raisebox{5pt}{\lefteqn{\scriptscriptstyle #1}}
\raisebox{-1pt}{$\scriptscriptstyle #2$}\hspace{-1pt}\right\}
}}

\def\minitbl#1#2{{
\left\{\hspace{-1pt}
\raisebox{4pt}{\lefteqn{\scriptscriptstyle #1}}
\raisebox{-2pt}{$\scriptscriptstyle #2$}\hspace{-1pt}\right\}
}}

\newcounter{theorem}
\newtheorem{pred}[theorem]{Proposition}
\newtheorem{zam}[theorem]{Remark}
\newtheorem{lem}[theorem]{Lemma}
\newtheorem{opred}[theorem]{Definition}
\newtheorem{teor}[theorem]{Theorem}

\begin{document}

\title{\parbox{1.0\hsize}{
\begin{flushright}
\normalsize\it
To L.D. Faddeev\\
with the respect and appreciation
\end{flushright}
\begin{center}
\huge
Cayley-Hamilton theorem for quantum matrix\\
algebras of $GL(m|n)$ type
\end{center}
 } }

\author{D.I.~Gurevich\thanks{
Dimitri.Gourevitch@univ-valenciennes.fr}\\
{\small\it
ISTV, Universit\'e de Valenciennes, 59304 Valenciennes,
France}\\ 
P.N.~Pyatov\thanks{pyatov@thsun1.jinr.ru}\\
{\small\it Bogoliubov Laboratory of Theoretical Physics,
JINR, 141980 Dubna, Russia}\\
P.A.~Saponov\thanks{Pavel.Saponov@ihep.ru}\\
{\small\it Department of Theoretical Physics, IHEP,
142281 Protvino,  Russia}}
\maketitle

\begin{abstract}
The classical Cayley-Hamilton identities are generalized to
quantum matrix algebras of the $GL(m|n)$ type.
\end{abstract}

\section{Introduction}

The well known Cayley-Hamilton theorem of classical matrix analysis
claims that any square matrix $M$ with entries from a field $\Bbb K$
satisfies the identity
\be
\lb{odin}
\Delta(M) \equiv 0,
\ee
where $\Delta(x) = \det (M-xI)$ is the characteristic polynomial of
$M$. Provided $\Bbb K$ be algebraically closed, the coefficients of
this polynomial are the elementary symmetric functions in the
eigenvalues of $M$. Below for the sake of simplicity we assume
$\Bbb K$ to be the field $\Bbb C$ of complex numbers.

In \cite{NT,PS,GPS,IOPS,IOP1} the Cayley-Hamilton theorem was
successively generalized to the case of the quantum matrix algebras
of the $GL(m)$ type. Any such an algebra is constructed by means of
R-matrix representations of Hecke algebras of the $GL(n)$ type.

In the works \cite{KT} the Cayley-Hamilton identity was established
for matrix superalgebras. The aim of the present paper consists
in a further generalization of the results of \cite{KT} to the
case of quantum matrix superalgebras.

Let us formulate our main result --- the Cayley-Hamilton theorem
for the quantum matrix algebras of the $GL(m|n)$ type. All the
definitions and notations are explained in the main text below.

{\bf Theorem.} {\it Let $M=\|M_i^j\|_{i,j=1}^N$ be the matrix of the
generators of the quantum matrix algebra ${\cal M}(R,F)$ defined by
a pair of compatible, strictly skew-invertible R-matrices $R$ and
$F$, where $R$ is a $GL(m|n)$ type R-matrix (see sections
\ref{sec:mrf1} and \ref{sec:mrf2}). Then in ${\cal M}(R,F)$ the
following matrix identity holds
\be
\lb{i1}
\sum_{i=0}^{n+m}  M^{\overline{m+n-i}}\, C_i\,\equiv 0\, ,
\ee
where $M^{\bar i}$ is the $i$-th power of the quantum matrix (see
definition (\ref{m-k})) and the coefficients $C_i$ are linear
combinations of the Schur functions (see definition (\ref{f-sh})).

The Schur functions entering the expression of $C_i$ are homogeneous
polynomials of the $(mn+i)$ order in the generators $M_i^j$.
Therefore, any matrix element in the left hand side of (\ref{i1}) is
a homogeneous polynomial of the $(mn + m +n)$ order in $M_i^j$.

The expression of the coefficient $C_i$ in terms of the Schur
functions can be graphically presented as follows
\be
\lb{i2}
C_i =
\sum_{k=\max\{0,(i-n)\}}^{\min\{i,m\}}
(-1)^k q^{(2k-i)}\;
\lefteqn{
\hspace{9mm}\mbox{\raisebox{14.7mm}{\tiny $n$ nodes}}}
\lefteqn{
\hspace{6mm}\overbrace{\hspace{16mm}\rule{0pt}{11.8mm}}}
\lefteqn{
\hspace{7mm}\mbox{\raisebox{-14mm}{\tiny $(i-k)$ nodes}}}
\lefteqn{
\hspace{5.5mm}\mbox{\raisebox{-10mm}{$\underbrace{
\hspace{12mm}}$}}}
\raisebox{3.2mm}{\tiny $m$}
\raisebox{2.5mm}{$\left\{\rule{0pt}{9mm} \right.$}
\hspace{-2pt}
\begin{tabular}{|c|cc|c|}\hline
&\multicolumn{2}{c|}{\hspace{-4mm}$\lefteqn\dots$}&\\ \hline
&&&\raisebox{-2pt}{$\smash{\vdots}$}\\ \cline{4-4}
\raisebox{3pt}{$\smash{\vdots}\!$}&
\multicolumn{2}{c|}{\hspace{-5mm}\raisebox{5pt}{$\smash{
\lefteqn\ddots}$}}
&\multicolumn{1}{c}{}\\[-12pt]
&\multicolumn{2}{c|}{}&\multicolumn{1}{c}{}\\ \cline{1-3}
&\multicolumn{1}{c|}{$\lefteqn\dots$ $\;\;$}&
\multicolumn{2}{c}{}\\ \cline{1-2}
\end{tabular}
\hspace{-1pt}
\raisebox{5.4mm}{$\left. \rule{0pt}{6mm} \right\}$}
\raisebox{5.7mm}{\tiny $k$}
\ee
where each Young diagram stands for the corresponding Schur
function.}
\medskip

If the matrices $F$ and $R$ are the transpositions in the superspace
$(m|n)$, then the corresponding algebra ${\cal M}(R,F)$ is a matrix
superalgebra. In this case identity (\ref{i1}) coincides with the
invariant Cayley-Hamilton relations obtained in \cite{KT}.
Note that in the supersymmetric case the linear size $N$ of the
matrix $M$ is connected with the parameters $m$ and $n$ by the
equation $N=m+n$. As was shown in \cite{G}, in general this
relation does not take place in the quantum case. Our proof of
identity (\ref{i1}) does not rely on such like connections.

The structure of the work is as follows. Section \ref{sec:h-alg}
contains definitions and necessary facts on the Hecke algebras.
In section \ref{sec:mrf} we define the quantum matrix algebra of the
$GL(m|n)$ type. Ibidem the notions of the characteristic subalgebra
and of the matrix power space are introduced. The last section is
devoted to the proof of the Cayley-Hamilton identity (\ref{i1}).

\medskip

The authors express their appreciation to Alexei Davydov, Dimitry
Leites, Alexander Molev and  Hovhannes Khudaverdian. Special
gratitude to our permanent co-authors Alexei Isaev and Oleg
Ogievetsky for numerous fruitful discussions and remarks. The work
of the second author (P.N.P.) was partially supported by the RFBR 
grant No. 03-01-00781.

\section{Some facts on Hecke algebras}
\label{sec:h-alg}

In this section we give a summary of facts on the structure of the
Hecke algebras of $A_{n-1}$ series. For their proofs the reader is
referred to the works \cite{DJ,M}. Detailed consideration of the
material and the list of references on the theme can also be found
in the reviews \cite{OP,R}. When considering the quantum matrix
algebras we shall use this section as a base of our technical tools.

\subsection{Young diagrams and tableaux}

First of all we recall some definitions from the partition theory.
Mainly we shall use the terminology and notations of the book
\cite{Mac}.

{\it A partition} $\lambda$ of a positive integer $n$ (notation:
$\lambda\vdash n$) is a non-increasing sequence of non-negative
integers $\lambda_i$ such that their sum is equal to $n$:
$$
\lambda = (\lambda_1,\lambda_2,\dots),\quad \lambda_i\ge
\lambda_{i+1}:\quad |\lambda|\equiv \sum\lambda_i = n.
$$
The number $|\lambda|=n$ is called {\it the weight} of the partition
$\lambda$. A convenient graphical image for a partition is 
{\it a Young diagram} (see \cite{Mac}) and below we shall identify
any partition with the corresponding Young diagram.

With each node of a Young diagram located at the intersection of the
$s$-th column and the $r$-th row we associate a number $c=s-r$ which
will be called {\it the content} of the node.

Define an inclusion relation $\subset$ on the set of all the Young
diagrams by stating that a diagram $\mu$ includes a diagram
$\lambda$ (symbolically $\lambda \subset \mu$) if after the
superposition the diagram $\lambda$ is entirely contained inside of
$\mu$. In other words, $\lambda \subset \mu \Leftrightarrow
\lambda_i\leq\mu_i\;  i =1,2,\dots$

Given a Young diagram $\lambda\vdash n$ one can construct $n!$
{\it Young tableaux} by writing into the nodes of the diagram all
integers from 1 to  $n$ in an arbitrary order. The tableau will be
noted by the symbol
\vspace{-2mm}
\be
\tbl{\lambda}{\alpha},\qquad  \alpha =1,2,\dots , n!,
\label{tabl}
\vspace{-2mm}
\ee
where the index $\alpha$ marks the different tableaux. The diagram
$\lambda$ will be called {\it the form} of a tableau
$\tbl{\lambda}{\alpha}$.

The set of all tableaux associated with a diagram $\lambda\vdash n$
can be naturally supplied with an action of the $n$-th order
symmetric group $S_n$. By definition, an element $\pi\in S_n$ acts
on a tableaux $\tbl{\lambda}{\alpha}$ by changing any number $i$ for
$\pi(i)$ in all the nodes of the tableau. The resulting tableau will
be denoted as $\tbl{\lambda}{\pi(\alpha)}$.

In the set of all the Young tableaux of the form $\lambda$ we shall
distinguish a subset of {\it standard} tableaux which satisfy an
additional requirement: the numbers written in a standard tableau
increases in any column from top to bottom and in any raw from left
to right. From now on we shall use only standard Young tableaux and
retain notation (\ref{tabl}) for them. However, the index $\alpha$
will
vary from 1 to an integer $d_\lambda$ which is equal to the total
number of all the standard Young tableaux of the same form
$\lambda$. Explicit expressions for $d_\lambda$ are given, for
example, in \cite{Mac}.

Though the subset of the standard tableaux is not closed with
respect to the action of the symmetric group, it nevertheless is a
cyclic set under such an action since an arbitrary standard Young
tableau of the form $\lambda\vdash n$ can be obtained from any other
standard tableau of the same form by the action of some element
$\pi$ from the group $S_n$ that is
$$
\forall\,\alpha,\beta\quad \exists\, \pi\in S_n:
\quad \tbl{\lambda}{\alpha}\mapsto\tbl{\lambda}{\beta}\equiv
\tbl{\lambda}{\pi(\alpha)}.
$$

The relation $\subset$ defined above for the Young diagrams can be
extended to the set of the standard tableaux. Namely, we say that
a tableau includes another one
\vspace{-2mm}
$$
\tbl{\lambda}{\alpha}\subset \tbl{\mu}{\beta},
\vspace{-3mm}
$$
if $\lambda\subset\mu$ and the numbers from $1$ to $|\lambda|$
occupy the same nodes in the tableau \tbl{\mu}{\beta} as they do in
\tbl{\lambda}{\alpha}.

\subsection{Hecke algebras: definition and the basis of matrix
units}

{\it The Hecke algebra} ${\cal H}_n(q)$ of $A_{n-1}$ series is a
unital associative ${\Bbb C}$-algebra generated by the set of
elements $\{\sigma_i\}_{i=1}^{n-1}$ subject to the relations
\ba
\lb{1}
\sigma_k \sigma_{k+1} \sigma_k = \sigma_{k+1} \sigma_k
\sigma_{k+1}\ , &\quad & k=1,\dots ,n-2 \; , \\
\lb{2}
\sigma_k \sigma_j = \sigma_j \sigma_k \; , &\quad & \mbox{if} \;\;
|k-j|\geq 2 \; , \\
\lb{3}
{\sigma_k}^2 = \mathbf{1} + (q-q^{-1}) \sigma_k\ , &\quad &
 k=1,\dots , n-1 \; ,
\ea
where $\mathbf{1}$ is the unit element. The last of the above
defining relations depends on a parameter
$q\in {\Bbb C}\setminus \{0\}$ and is referred to as {\it the Hecke
condition}.

${\cal H}_n(1)$ is identical with the group algebra ${\Bbb C}[S_n]$
of the symmetric group $S_n$. In case
\be
k_q := {q^k - q^{-k}\over q-q^{-1}}\neq 0\, ,\quad \mbox{for all}
\quad k=2,3,\dots ,n,
\label{q-num}
\ee
the algebra ${\cal H}_n(q)$ is semisimple and is isomorphic to
${\Bbb C}[S_n]$. In what follows we shall suppose condition
(\ref{q-num}) to  be fulfilled for all
$k\in {\Bbb N}\setminus\{0\}$. Under such an assumption the algebra
${\cal H}_n(q)$ is isomorphic to the direct sum of matrix algebras
\be
{\cal H}_n(q) \cong\bigoplus_{\lambda\vdash n}{\rm Mat}_{d_\lambda}
({\Bbb C}).
\label{iso}
\ee
Here ${\rm Mat}_m({\Bbb C})$ stands for an associative
${\Bbb C}$-algebra generated by $m^2$ matrix units
$E_{ab}$ $1\le a,b\le m$ which satisfy the following multiplication
law
$$
E_{ab}\,E_{cd} = \delta_{bc}\, E_{ad}.
$$

Due to isomorphism (\ref{iso}) one can construct the set of elements
$E^{\lambda}_{\alpha\beta}\in {\cal H}_n(q)$,
$\forall\,\lambda\vdash n$, $1\le \alpha,\beta\le d_\lambda$,
which are the images of matrix units generating the term
${\rm Mat}_{d_\lambda}$ in the direct sum (\ref{iso})
\be
E^\lambda_{\alpha\beta}\,E^\mu_{\gamma\tau} = \delta^{\lambda\mu}
\delta_{\beta\gamma}\,E^\lambda_{\alpha\tau}.
\label{me-alg}
\ee
The diagonal matrix units
$E^\lambda_{\alpha\alpha}\equiv E^\lambda_{\alpha}$ are the
primitive idempotents of the algebra ${\cal H}_n(q)$ and the set
of all $E^\lambda_{\alpha\beta}$ forms a linear basis in
${\cal H}_n(q)$.

Obviously, the matrix units $E^\lambda_{\alpha\beta}$ can be
realized in many different ways. The construction presented in
\cite{M,R} (see also \cite{OP}) satisfies a number of useful
additional relations. It is this type of the basis of matrix
units which we shall use below. Now we shortly outline its main
properties.

Firstly, in this basis there are known remarkably simple expressions
for the elements which transform different diagonal matrix units
into each other. In order to describe these elements we need a few
new definitions.

With each generator $\sigma_i$ of the algebra $H_n(q)$ we associate
the element
\be
\sigma_i(x):=\sigma_i+\frac{q^{-x}}{x_q}\,\mathbf{1}, \quad
1\le i\le n-1.
\label{yb-form}
\ee
Here $x$ is an arbitrary non-negative integer, the symbol $x_q$
being defined in (\ref{q-num}). The elements $\sigma_i(x)$ obey
the relations
\ba
\sigma_i(x)\sigma_{i+1}(x+y)\sigma_i(y) &=&
\sigma_{i+1}(y)\sigma_{i}(x+y)\sigma_{i+1}(x),\qquad 1\le i\le n-2,
\label{ybe-par}
\\[2mm]
\lb{unitarity}
\sigma_i(x)\sigma_i(-x) &=&
{(x+1)_q (x-1)_q\over x_q^2}\, {\bf 1},\qquad\qquad\;
1\le i\le n-1.
\ea

Consider an arbitrary standard Young tableau \tbl{\lambda}{\alpha}
corresponding to a partition $\lambda\vdash n$. For each integer
$1\le k\le n-1$ we define the number $\mtbl{\lambda}{\alpha}{k}$
to be equal to the difference of the contents of nodes containing
the numbers $k$ and $(k+1)$
\be
\lb{content}
\mtbl{\lambda}{\alpha}{k}:= c(k) - c(k+1).
\ee
Let a tableau \tbl{\lambda}{\beta} (not obligatory standard) differs
from the initial tableau \tbl{\lambda}{\alpha} only in the
transposition of the nodes with the numbers $k$ and $(k+1)$
\be
\lb{beta}
\beta = \pi_k(\alpha),\quad 1\le k\le n-1,
\ee
$\pi_k\in S_n$ being the transposition of $k$ and $(k+1)$. If the
tableau \tbl{\lambda}{\pi_k(\alpha)} is standard then the following
relations hold
\ba
\lb{me-con1}
\sigma_k(\mtbl{\lambda}{\alpha}{k})\,E^\lambda_{\alpha} &=&
E^\lambda_{\pi_k(\alpha)}\,\sigma_k(-\mtbl{\lambda}{\alpha}{k})\,,
\\
\lb{me-con1a}
E^\lambda_{\alpha}\,\sigma_k(\mtbl{\lambda}{\alpha}{k})&=&
\sigma_k(-\mtbl{\lambda}{\alpha}{k})\,E^\lambda_{\pi_k(\alpha)}\,.
\ea
Now, take into account the cyclic property of the set of standard
Young tableaux under the action of the symmetric group and the fact
that an arbitrary permutation can be (non-uniquely) expanded into a
product of transpositions. This allows one to connect any pair
$E^\lambda_\alpha$ and $E^\lambda_\beta$ by a chain of
transformations similar to (\ref{me-con1}), (\ref{me-con1a}).
Moreover, relations (\ref{ybe-par}), (\ref{unitarity}) ensure the
compatibility of all possible formulae expressing dependencies
between $E^\lambda_\alpha$ and $E^\lambda_\beta$.

If the tableau \tbl{\lambda}{\pi_k(\alpha)} is non-standard then
\be
\sigma_k(\mtbl{\lambda}{\alpha}{k})\,E^\lambda_{\alpha} =
E^\lambda_{\alpha}\,\sigma_k(\mtbl{\lambda}{\alpha}{k}) = 0.
\label{me-con2}
\ee

The expressions for non-diagonal matrix units follow from
(\ref{me-con1}), (\ref{me-con1a}). Assuming both the tableaux
\tbl{\lambda}{\alpha} and \tbl{\lambda}{\pi_k(\alpha)} to be
standard, one gets
\ba
\lb{me-con3}
E^\lambda_{\alpha\pi_k(\alpha)} &:=&
\omega(\mtbl{\lambda}{\alpha}{k})\,
E^\lambda_{\alpha}\,\sigma_k(\mtbl{\lambda}{\alpha}{k})\, ,
\\
\lb{me-con4}
E^\lambda_{\pi_k(\alpha)\alpha} &:=&
\omega(\mtbl{\lambda}{\alpha}{k})\,
\sigma_k(\mtbl{\lambda}{\alpha}{k})\,E^\lambda_{\alpha}\, ,
\ea
where the normalizing coefficient $\omega(\ell)$ satisfies the
relation
\be
\lb{w-l}
\omega(\ell)\, \omega(-\ell)=
\frac{\ell_q^2}{(\ell+1)_q(\ell-1)_q}\ .
\ee
In particular, one can set $\omega(\ell):=\ell_q/(\ell_q+1)$.
Note that in case the tableau \tbl{\lambda}{\pi_k(\alpha)} is a
standard one, we have $\mtbl{\lambda}{\alpha}{k}\not=\pm1$.
Therefore, the above expression for the normalizing coefficient
is always correctly defined.
\medskip

To formulate the second property of the matrix units which will be
used below let us consider the chain of Hecke algebra embeddings
\be
{\cal H}_2(q)\hookrightarrow\dots\hookrightarrow
{\cal H}_k(q)\hookrightarrow
{\cal H}_{k+1}(q)\hookrightarrow\dots ,
\lb{h-emb}
\ee
which are defined on generators by the formula
\be
\lb{h-emb2}
{\cal H}_{k}(q)\ni \sigma_i \mapsto \sigma_{i+1}\in {\cal H}_{k+1}(q),
\quad  i=1,\dots ,k-1.
\ee
In each subalgebra ${\cal H}_k(q)$ one can choose the basis of
the matrix units $E^\lambda_{\alpha\beta}$ in such a way that the
diagonal matrix units $E_\alpha^\lambda\in {\cal H}_k(q)$ can be
decomposed into the sum of diagonal matrix units belonging to any
enveloping algebra ${\cal H}_m(q)\supset {\cal H}_k(q)$
\be
\lb{h-emb3}
E^{\lambda\vdash k}_{\alpha} =
\sum_{\minitbl{\lambda}{\alpha}\subset\minitbl{\mu}{\beta}}
E^{\mu\vdash m}_\beta\, , \quad m\ge k.
\ee
Here the summation goes over all the standard Young tableaux
\tbl{\mu}{\beta}, $\mu\vdash m$, containing the standard tableau
\tbl{\lambda}{\alpha}.

\section{Quantum matrix algebra ${\cal M}(R,F)$}
\label{sec:mrf}

\subsection{R-matrix representations of the Hecke algebra}
\label{sec:mrf1}

Let $V$ be a finite dimensional vector space over the field of
complex numbers, $\dim V = N$. With any element
$X\in {\rm End}(V^{\otimes p})$, $p= 1,2,\dots ,$ we associate
a sequence of endomorphisms $X_i\in {\rm End}(V^{\otimes k})$,
$k= p,p+1,\dots ,$ according to the rule
\be
\lb{endo}
X_i = {\rm Id}_V^{i-1}\otimes X\otimes {\rm Id}_V^{k-p-i+1},
\quad 1\leq i\leq k-p+1,
\ee
where ${\rm Id}_V$ is the identical automorphism of $V$.

An operator $R\in {\rm Aut}(V^{\otimes 2})$ satisfying the
{\em Yang-Baxter equation}
\be
R_1R_2R_1 = R_2R_1R_2\,,
\lb{ybe}
\ee
will be called {\em the R-matrix}. One can easily see that any
R-matrix obeying the condition
\be
(R-q\,{\rm Id}_{V^{\otimes 2}})(R+q^{-1}\,{\rm Id}_{V^{\otimes 2}})
= 0, \quad q\in {\Bbb C}\setminus \{0\}
\lb{h-cond}
\ee
generates representations $\rho_R$ for the series of the Hecke
algebras ${\cal H}_k(q)$, $k=2,3,\dots$
\be
\rho_R: {\cal H}_k(q)\rightarrow {\rm End}(V^{\otimes k})\qquad
\sigma_i\mapsto \rho_R(\sigma_i) = R_i\qquad 1\le i\le k-1,
\lb{h-rep}
\ee
where $\sigma_i$ are the generators of ${\cal H}_k(q)$ (\ref{1}) --
(\ref{3}). Such like R-matrices will be called {\em the Hecke type}
(or simply {\em Hecke}) R-matrices and the corresponding
representations $\rho_R$  --- {\em the R-matrix} ones.

Let $R$ be a Hecke R-matrix. Suppose that the corresponding
representations
$\rho_R: {\cal H}_k(q)\rightarrow \mbox{End}(V^{\otimes k})$ are
faithful for all $2\le k <(m+1)(n+1)$ while the representation
$\rho_R: {\cal H}_{(m+1)(n+1)}(q)\rightarrow \mbox{End}(V^{\otimes
(m+1) (n+1)})$ possesses a kernel, generated by any matrix unit
corresponding to the rectangular Young diagram $((n+1)^{(m+1)})$
that is
\be
\begin{array}{l}
\rho_R(E^{((n+1)^{(m+1)})}_\alpha) = 0\\
\rule{0pt}{6mm}
\rho_R(E^{\mu}_\alpha)\neq 0, \qquad
\mbox{for any}\; \mu\vdash (m+1)(n+1),
\quad \mu\neq ((n+1)^{(m+1)}).
\end{array}
\lb{glmn-rep}
\ee
Such an $R$-matrix will be referred to as an R-matrix of
{\em the $GL(m|n)$ type}.

\subsection{R-trace and pairs of compatible R-matrices}
\label{sec:mrf2}

Let $\{v_i\}_{i=1}^N$ be a fixed basis in the space $V$. In the
basis $\{v_i\otimes v_j\}_{i,j = 1}^N$ an operator
$X\in {\rm End}(V^{\otimes 2})$ is defined by its matrix
$X_{ij}^{kl}$: $X(v_i\otimes v_j):=\rule{0pt}{4mm}
\sum_{k,l=1}^N X_{ij}^{kl}\, v_k\otimes v_l$.

The operator $X\in {\rm End}(V^{\otimes 2})$ will be called
{\it skew-invertible} if there exists a matrix ${\Psi^X}^{kl}_{ij}$
such that
\be
\sum_{a,b = 1}^N X_{ia}^{kb} {\Psi^X}^{al}_{bj} = \delta^l_i\,
\delta^k_j.
\lb{s-inv}
\ee
One can easily see that the property of skew-invertibility is
invariant with respect to the change of the basis $\{v_i\}$ and,
therefore, to the matrix ${\Psi^X}^{kl}_{ij}$ there corresponds an
operator $\Psi^X\in {\rm End (V^{\otimes 2})}$. Introduce a notation
\be
D^X:={\rm Tr}_{(2)}\Psi^X\in \mbox{End}(V),
\ee
where the symbol ${\rm Tr}_{(2)}$ means taking the trace over the
second space in the product \mbox{$V\!\otimes\! V$}. A
skew-invertible operator $X$ will be called {\em strictly}
skew-invertible in case the operator $D^X$ is invertible.

Consider the set of $N\times N$ matrices ${\rm Mat}_N(W)$, their
entries being vectors of a linear space $W$ over $\Bbb C$. The set
${\rm Mat}_N(W)$ is endowed with the structure of a linear space
over $\Bbb C$ in the obvious way.

Let $R$ be a skew-invertible R-matrix. The linear map
$$
{\rm Tr\str{-1.3}}_R:\; {\rm Mat}_N(W)\,\rightarrow \, W ,
$$
defined by
\be
{\rm Tr\str{-1.3}}_R(M) = \sum_{i,j=1}^N{D^R}_i^jM_j^i\qquad
M\in{\rm Mat}_N(W)\,,
\lb{r-sled}
\ee
will be called the calculation of {\em R-trace}.

\begin{zam}{\rm If as an R-matrix one chooses the transposition
operator in the superspace of $(m|n)$ type (see definition
(\ref{superP})), then the operation ${\rm Tr\str{-1.3}}_R$
coincides with the calculation of the supertrace.}
\end{zam}

An ordered pair $\{R, F\}$ of the two R-matrices $R$ and $F$ will
be called {\em a compatible pair} if the following relations are
satisfied
\be
\begin{array}{c}
R_1F_2F_1 = F_2F_1R_2\\
R_2F_1F_2 = F_1F_2R_1.
\end{array}
\lb{sovm}
\ee

Now we list some properties of the compatible pairs of R-matrices.
\begin{enumerate}
\item If $R$ is a skew-invertible R-matrix and $\{R, F\}$ is a
compatible pair, then for all $M\in{\rm Mat}_N(W)$ the following
relation holds
\be
\Tr{2}(F_1^{\pm 1}M_1F_1^{\mp 1}) = {\rm Id}_V\,{\rm Tr\str{-1.3}}_R
M \quad \mbox{where} \quad M_1 = M\otimes {\rm Id}_V.
\lb{tr2-1}
\ee

\item
In case $R$ is skew-invertible, then
\be
[ R_1, D_1^R D_2^R ] = 0\, .
\lb{rdd}
\ee
A direct consequence of this formula is a {\em cyclic property} of
the R-trace
\be
\Tr{12}(R_1U) = \Tr{12}(UR_1),\quad U\in {\rm Mat}_N(W)^{\otimes 2}.
\lb{cikl}
\ee

\item If $\{R, F\}$ is a compatible pair, then $R_f:=F^{-1}R^{-1}F$
is an R-matrix and the pair $\{R_f, F\}$ is compatible.

\item Let $\{R, F\}$ be a compatible pair, R-matrix $R_f$ is
skew-invertible and matrices $R$ and $F$ are strictly
skew-invertible. Then the map $\phi: {\rm Mat}_N(W) \rightarrow
{\rm Mat}_N(W)$ defined by 
\be
\phi(M) := \Tr{2}(F_1M_1F_1^{-1}R_1), \quad 
M\in {\rm Mat}_N(W),
\lb{fi}
\ee
is invertible. An explicit form of its inverse $\phi^{-1}$
is as follows (see \cite{OP2})
\be
\phi^{-1}(M) = {\rm Tr \str{-1.3}}_{{R_f}^{\mbox{
\scriptsize $(2)$}}} (F_1^{-1}M_1R_1^{-1}F_1).
\lb{fi-inv}
\ee
\end{enumerate}

\subsection{Quantum matrix algebras: definition}

Consider a linear space ${\rm Mat}_N(W)$ and define a series of
linear mappings ${\rm Mat}_N(W)^{\otimes k}\rightarrow
{\rm Mat}_N(W)^{\otimes (k+1)}$, $k\ge 1$, by the following
recurrent relations
\be
M_{\overline 1}:=M, \quad M_{\overline{k}} \mapsto
M_{\overline{k+1}}:= F_kM_{\overline{k}}F_k^{-1}, \quad
M_{\overline k}\in {\rm Mat}_N(W)^{\otimes k},
\lb{kopii}
\ee
where $M$ is an arbitrary matrix from ${\rm Mat}_N(W)$ and $F$ is a
fixed element of ${\rm Aut}(V\otimes V)$.

As the linear space $W$ we choose an associative algebra ${\cal A}$
over $\Bbb C$ freely generated by the unit element and $N^2$
generators $M_i^j$
$$
{\cal A} = {\Bbb C}\langle M_i^j\rangle\quad 1\le i,j\le N.
$$

\begin{opred}
\lb{opredelenie}
Let $\{R,F\}$ be a compatible pair of strictly skew-invertible
R-matrices, the R-matrix $R_f=F^{-1}R^{-1}F$ being skew-invertible.
By definition, the quantum matrix algebra ${\cal M}(R,F)$ is a
quotient of the algebra ${\cal A}$ over a two-sided ideal generated
by entries of the following matrix relation
\be
R_1M_{\overline 1}M_{\overline 2} - M_{\overline 1}
M_{\overline 2}R_1 = 0,
\label{qma}
\ee
where $M = \|M_i^j\|\in {\rm Mat}_N({\cal A})$ and the matrices
$M_{\overline k}$ are constructed according to (\ref{kopii}) with
the R-matrix $F$.
\end{opred}

\begin{zam}{\rm The definition of the quantum matrix algebra
presented here differs from that given in \cite{IOP1}. Nevertheless,
under the additional assumptions on $R$ and $F$ made above, both
definitions are equivalent. As was shown in \cite{IOP2}, the matrix
$\widetilde M = (D^F)^{-1}MD^R$ obeys relations imposed on the
generators of the quantum matrix algebra in the work \cite{IOP1}.}
\end{zam}

\begin{lem}
The matrix $M$ composed of the generators of the quantum matrix
algebra ${\cal M}(R,F)$ satisfies the relations
\be
R_k\,M_{\overline k}\,M_{\overline {k+1}} = M_{\overline k}\,
M_{\overline {k+1}}\, R_k,
\lb{rmm-k}
\ee
where the matrices $M_{\overline k}$, $M_{\overline{k+1}}$ are
defined by rules (\ref{kopii}) with the R-matrix $F$.
\end{lem}

The proof of the lemma is presented in \cite{IOP1}.

\subsection{Characteristic subalgebra and Schur functions}

Further in this section we consider the quantum matrix algebras
${\cal M}(R,F)$ defined by a Hecke R-matrix $R$. We shall refer to
them as the matrix algebras of {\em the Hecke type}.

Let a subspace ${\rm Char}(R,F)\subset{\cal M}(R,F)$ be a linear
span of the unit and of the following elements
\be
\lb{char}
y(x^{(k)}) = \Tr{1\dots k}(M_{\overline 1}\dots M_{\overline
k}\,\rho_R(x^{(k)}))\quad k =1,2,\dots ,
\ee
where $x^{(k)}$ runs over all elements of ${\cal H}_k(q)$.
The symbol $\Tr{1\dots k}$ stands for the calculation of the R-trace
over the spaces from the first to the $k$-th ones.

\begin{pred}
Let ${\cal M}(R,F)$ be a quantum matrix algebra of the Hecke type.
The subspace ${\rm Char}(R,F)$ is a commutative subalgebra of
${\cal M}(R,F)$. The subalgebra ${\rm Char}(R,F)$ will be called
the characteristic subalgebra.
\end{pred}
The detailed proof is given in  \cite{IOP1}. It is based, in
particular, on the following technical result.
\begin{lem}
Consider an arbitrary element $x^{(k)}$ of the algebra
${\cal H}_k(q)$. Let $x^{(k)\uparrow i}$ denotes the image of
$x^{(k)}$ in the algebra ${\cal H}_{k+i}(q)$ under the embedding
${\cal H}_k(q) \hookrightarrow {\cal H}_{k+i}(q)$ defined by
(\ref{h-emb2}). Let $\{R,F\}$ be a compatible pair of R-matrices,
the Hecke R-matrix $R$ being skew-invertible. Then the following
relation holds true
\be
\lb{char1}
\Tr{i+1,\dots ,i+k}(M_{\overline{i+1}}\dots
M_{\overline{i+k}}\ \rho_R(x^{(k)\uparrow i}))\, =\,
{\rm Id}_{V^{\otimes i}}\  y(x^{(k)}).
\ee
\end{lem}
Formula (\ref{char1}) will be repeatedly used in the sequel. It is
immediate consequence of the property (\ref{tr2-1}) of the
compatible pairs $\{R,F\}$ and its proof is also presented in
\cite{IOP1}.
\medskip

Take into consideration two sets $\{p_k(M)\}$ and $\{s_\lambda(M)\}$
of elements of the characteristic subalgebra. The elements of these
sets will be referred to as {\em the power sums} and {\em the Schur
functions} respectively. For an arbitrary integer $k\geq 1$ the
power sums are defined by
\be
p_k(M):= \Tr{1\dots k}(M_{\overline 1}\dots
M_{\overline k}\,R_{k-1}\dots R_1).
\lb{st-sm}
\ee
Given a Hecke type R-matrix $R$ and an arbitrary partition
$\lambda\vdash k$, $k=1,2,\dots$, we define the Schur functions
as follows
\be
s_0(M) := 1,\quad
s_\lambda(M) := \Tr{1\dots k}(M_{\overline 1}\dots
M_{\overline k}\, \rho_R(E^\lambda_\alpha)).
\lb{f-sh}
\ee
The right hand side of (\ref{f-sh}) does not actually depend on the
index $\alpha$ of the matrix unit. Indeed, for all
$x^{(k)}\in {\cal H}_k(q)$ the matrices $\rho_R(x^{(k)})$ and
$(M_{\overline 1}\dots M_{\overline k})$ commute in (\ref{char})
due to (\ref{rmm-k}). Then, taking into account the cyclic property
of the R-trace, one gets
\be
\lb{xox}
y\left(\sigma_k(\ell)\, E^\lambda_\alpha\,
\sigma^{-1}_k(\ell)\right)\ =\  y(E^\lambda_\alpha)\, , \quad
y(E^\lambda_{\pi_k(\alpha)\alpha})= w(\ell)\,
y\left(E^\lambda_{\pi_k(\alpha)}\,\sigma_k(\ell)\,
E^\lambda_\alpha\right)= 0,
\ee
where $\ell := \mtbl{\lambda}{\alpha}{k}$. Next, with the use of
(\ref{yb-form}), (\ref{unitarity}), (\ref{me-con4}) and (\ref{w-l})
relations (\ref{me-con1}) can be rewritten in the form
\be
\lb{opsa}
E^\lambda_{\pi_k(\alpha)}\ =\ \sigma_k(\ell)\, E^\lambda_\alpha\,
\sigma^{-1}_k(\ell)\ -\ (q-q^{-1})\, w(-\ell)\,
E^\lambda_{\pi_k(\alpha)\alpha}\,.
\ee
This leads to $y(E^\lambda_{\pi_k(\alpha)}) = y(E^\lambda_\alpha)$.
At last, the fact that the set of the standard Young tableaux is
cyclic with respect to the action  of the symmetric group guarantees
that definition (\ref{f-sh}) is non-contradictory.

\begin{pred}
\lb{pred:basis}
$\!\!\!$ For any quantum matrix algebra ${\cal M}(\!R,\!F)$ of the
Hecke type
\begin{itemize}
\item[]{i)}
Its characteristic subalgebra ${\rm Char}(R,F)$ is generated by the
unit and by the power sums $p_k(M)$, $k=1,2,\dots$;
\item[]{ii)}
The set of Schur functions $s_\lambda(M)$, $\lambda\vdash k$,
$k=0,1,\dots$, forms a linear basis in ${\rm Char}(R,F)$.
\end{itemize}
\end{pred}

\noindent
{\bf Proof.}~The first part of the proposition is proved in
\cite{IOP1}. To prove the second part, we note that, basing on the
arguments given to justify relations (\ref{f-sh}), one can get
the following equalities
\be
\lb{nul}
y(E^\lambda_{\alpha\beta})\ =\ \delta_{\alpha\beta}\,
y(E^\lambda_{\alpha}) ,\quad y(E^\lambda_{\alpha})\ =\
y(E^\lambda_{\beta}) , \qquad \forall\;\alpha , \beta .
\ee
Now the statement on a linear basis of the Schur functions follows
from the fact that the matrix units
$\{E^\lambda_{\alpha\beta}\,|\, \lambda\vdash k\}$ form the linear
basis in the algebra ${\cal H}_k(q)$.\hfill
\raisebox{3pt}{\framebox[3mm]{\rule{0pt}{1mm}}}

\subsection{Powers of quantum matrices}

Let a subspace ${\rm Pow}(R,F)\subset {\rm Mat}_N({\cal M}(R,F))$
be the linear span of the unit matrix and of the following matrices
\be
M^{(x^{(k)})}:= \Tr{2\dots k}(M_{\overline 1}\dots M_{\overline k}\,
\rho_R(x^{(k)})) \qquad  k=1,2,\dots ,
\lb{x-step}
\ee
where $x^{(k)}$ runs over all elements of ${\cal H}_k(q)$. The
matrix $M^{(x^{(k)})}$ will be called {\em the $x^{(k)}$-th power}
of the matrix $M$, and the space ${\rm Pow}(R,F)$ will be referred
to as {\em the space of matrix powers}.
\begin{pred}
\lb{pred:8}
${\rm Pow}(R,F)$ is a right module over the characteristic
subalgebra ${\rm Char}(R,F)$.
\end{pred}
With the use of (\ref{char1}) this proposition can be proved
analogously to the proof of the commutativity of ${\rm Char}(R,F)$
(see \cite{IOP1}).
\medskip

In the space ${\rm Pow}(R,F)$ we choose a special set of matrices
constructed by the following rule. Having fixed a standard Young
tableau \tbl{\lambda}{\alpha} of the form $\lambda\vdash k$,
$k=1,2,\dots$, we define the matrices
\be
M^{((1);1)}:=M,\quad
M^{(\lambda;i)}:=\Tr{2\dots k}(M_{\overline 1}\dots
M_{\overline k}\,
\rho_R(E^\lambda_\alpha)), \quad k=2,3,\dots ,
\lb{m-lam}
\ee
where the index $i$ in $M^{(\lambda;i)}$ is the number of the raw of
the tableau \tbl{\lambda}{\alpha} (when counting from the top) which
contains $k$ --- the largest integer among those filling the
tableau. Take, as an example, the standard tableau
$$
\hspace{32mm}
\begin{tabular}{|c|c|c|c|}\hline
$\smash{\scriptstyle 1}$&$\smash{\scriptstyle 3}$&
$\smash{\scriptstyle 4}$&$\smash{\scriptstyle 6}$\\ \hline
$\smash{\scriptstyle 2}$&$\smash{\scriptstyle 7}$&
\multicolumn{2}{c}{}\\ \cline{1-2}
$\smash{\scriptstyle 5}$&\multicolumn{3}{c}{}\\ \cline{1-1}
\end{tabular}.
$$
Since the largest integer $k=7$ is contained in the second raw
then $i=2$ for the diagram involved.

The matrix $M^{(\lambda;i)}$ defined in (\ref{m-lam}) does not
depend on the positions of the other integers which are less than
$k$. The proof is analogous to that used when arguing in favour of
the unambiguity of the Schur function definition (see the paragraph
next to the formula (\ref{f-sh})). One should only take into account
an additional circumstance. If two matrix units $E^\lambda_\alpha$
and $E^\lambda_\beta$ $\lambda\vdash k$ have the number $k$ standing
in the same nodes of the corresponding tableaux, then the chain of
transformations linking these matrix units contains only the elements
$x^{(k-1)}$ of the subalgebra ${\cal H}_{k-1}(q)$. Just with respect
to the images $\rho_R(x^{(k-1)})\in \mbox{Id}_V\otimes
\mbox{Aut}(V^{\otimes (k-1)})$ the operation $\Tr{2\dots k}$
possesses the cyclic property.

Besides $M^{(\lambda;i)}$, consider also a series of matrices
$M^{\overline k}$:
\be
M^{\overline 1}:= M,\quad M^{\overline k}:= 
\Tr{2\dots k}(M_{\overline 1}\dots M_{\overline k}
\,R_{k-1} \dots R_1), \quad k=2,3,\dots .
\lb{m-k}
\ee
The matrices $M^{(\lambda;i)}$ and $M^{\overline k}$ will be called
respectively {\em the $(\lambda;i)$-th} and {\em the $k$-th powers}
of the matrix $M$.

\begin{pred}
\lb{2bazisa}
For any quantum matrix algebra ${\cal M}(R,F)$ of the Hecke type:
\begin{itemize}
\item[]{i)}
The unit matrix and the family of matrices $M^{\overline k}$, 
$k=1,2,\dots$ is a generating set of the right
${\rm Char}(R,F)$-module ${\rm Pow}(R,F)$;
\item[]{ii)}
The unit matrix and the family of matrices $M^{(\lambda;i)}$,
$\lambda\vdash k$, $k=1,2,\dots$ where $i$ runs over all the possible
values for each partition $\lambda\vdash k$ is a basis of the linear
space ${\rm Pow}(R,F)$.
\end{itemize}
\end{pred}

The proof of items {\it i)} and {\it ii)} of the above proposition
can be given by a quite obvious generalization of the proof of the
corresponding items of proposition \ref{pred:basis}.
\medskip

To end the section, we present another form of the $k$-th power of
the matrix $M$. Basing on the property (\ref{tr2-1}) of compatible
pairs one can show that definition (\ref{m-k}) is equivalent to the
following iterative formulae
\be
\lb{m-k2}
M^{\overline{0}} = {\rm Id}_V, \quad
M^{\overline{k}} = M\cdot\phi(M^{\overline{k-1}}),
\ee
where the map $\phi$ is given by (\ref{fi}) and the dot stands for
the usual matrix product. Relations (\ref{m-k2}) can be interpreted
as the result of the consecutive application of the operator
${\sf M}: {\rm Mat}_N({\cal M}(R,F))\rightarrow
{\rm Mat}_N({\cal M}(R,F))$
$$
{\sf M}(X):=M\cdot \phi(X), \quad  X\in{\rm Mat}_N({\cal M}
(R,F))
$$
to the unit matrix ${\rm Id}_V$. Indeed, since
$\phi({\rm Id}_V) = {\rm Id}_V$ one gets
$$
{\sf M}({\rm Id}_V) = M,\quad M^{\overline k} = {\sf M}(
M^{\overline {k-1}}) =\dots ={\sf M}^k({\rm Id}_V).
$$

\subsection{Quantum matrix algebras: examples}

Let $P:V^{\otimes 2}\rightarrow V^{\otimes 2}$ be the transposition
operator:
$$
P(v_1\otimes v_2) = v_2\otimes v_1
$$
and R-matrix $R$ be strictly skew-invertible. The pair $\{R,P\}$ is
compatible and it defines a quantum matrix algebra ${\cal M}(R,P)$.
Denote the matrix of its generators by the letter $T$. In the case
involved, $T_{\overline k} = T_k$ (see definitions (\ref{endo}) and
(\ref{kopii})) and commutation relations (\ref{qma}) take the form
\be
R_1T_1T_2 - T_1T_2R_1 = 0.
\lb{RTT}
\ee
Suppose additionally that $R$ is of the $GL(m)$ type. The standard
example of such an R-matrix is the Drinfeld-Jimbo R-matrix obtained
by the quantization of classical groups of the $A_m$ series. The
corresponding algebra ${\cal M}(R,P) = {\rm Fun}_q(GL(N))$ is the
quantization of the algebra of functions on the general linear group
\cite{FRT}. Note that in general $m\not=N$. Examples of
$N^2\times N^2$ R-matrices of the $GL(m)$ type such that $m\not =N$
are constructed in \cite{G}.

In case $R$ is an R-matrix of the $GL(m)$ type (including the
Drinfeld-Jimbo R-matrix) the matrix $T$ of generators of
${\cal M}(R,P)$ obeys a polynomial Cayley-Hamilton identity of the
$m$-th order \cite{IOPS,IOP1}. The coefficients of the polynomial
are proportional to the Schur functions $s_{(1^k)}(T)$,
$0\le k\le m$. At the classical limit $q\to 1$ the Drinfeld-Jimbo
R-matrix tends to the transposition matrix $\lim_{q\to 1 }R = P$ and
algebra (\ref{RTT}) becomes the commutative matrix algebra
${\cal M}(P,P)$. The corresponding Cayley-Hamilton identity
coincides with (\ref{odin}). In particular, the Schur functions
$s_{(1^k)}(T)$ turn into the elementary symmetric functions in the
eigenvalues of the matrix $T$.

The other example of a compatible pair is given by $\{R,R\}$ where
$R$ is a strictly skew-invertible R-matrix. The commutation
relations (\ref{qma}) for the generators of the quantum matrix
algebra ${\cal M}(R,R)$ read
\be
R_1L_1R_1L_1  - L_1R_1L_1R_1 = 0,
\lb{REA}
\ee
where $L_1 = L\otimes {\rm Id}_V$ and $L = \|L_i^j\|$ is the matrix
composed of the algebra generators. The algebra ${\cal M}(R,R)$ is
called {\em the reflection equation algebra} \cite{KS}. Having first
appeared in the theory of integrable systems with boundaries, the
algebra has then found applications in the differential geometry of
the quantum groups (see, for example, \cite{IP,FP,GS,GS2}).

In this case the power of the quantum matrix coincides with the
usual matrix power $L^{\overline k} = L^k$. The characteristic
subalgebra ${\rm Char}(R,R)$ belongs to the center of
${\cal M}(R,R)$. The Cayley-Hamilton identity for the reflection
equation algebra with Drinfeld-Jimbo R-matrix was obtained in
\cite{NT,PS}. Its generalization to the case of an arbitrary
R-matrix of the $GL(m)$ type was found in \cite{GPS}.

Supposing $R$ to be a strictly skew-invertible R-matrix of the Hecke
type, let us make a linear shift of generators of ${\cal M}(R,R)$
$$
L\mapsto {\cal L} = L + \frac{1}{q-q^{-1}}{\rm Id}_V .
$$
After such a change relations (\ref{REA}) transforms into the
quadratic-linear form
\be
R_1{\cal L}_1R_1{\cal L}_1 - {\cal L}_1R_1{\cal L}_1R_1
= R_1{\cal L}_1 - {\cal L}_1 R_1.
\lb{mREA}
\ee
The above basis of generators is convenient in constructing the
quantum analogs of orbits of the coadjoint action of a Lie group
on the space dual to its Lie algebra. The Cayley-Hamilton identity
for this basis was obtained in \cite{GS}. If $R$ is the
Drinfeld-Jimbo R-matrix, then in the limit $q\to 1$ commutation
relations (\ref{mREA}) of the ${\cal M}(R,R)$ generators pass to
that of the generators of the universal enveloping algebra $U(gl_N)$
and the Cayley-Hamilton identity transforms into the well known
identity for the matrix composed of the $U(gl_N)$ generators (see,
for example, \cite{Gou}).

Let $P_{m|n}$ be the transposition operator on the superspace of the
$(m|n)$ type
\be
\lb{superP}
P_{m|n}(v_1\otimes v_2) = (-1)^{|v_1||v_2|} v_2\otimes v_1,
\ee
where $v_i$ is a homogeneous element of the superspace and $|v_i|$
denotes the corresponding parity. The algebra ${\cal M}(R,R)$ with
$R=P_{m|n}$ was considered in \cite{KT}. Note that the matrix
superalgebra ${\cal M}(P_{m|n},P_{m|n})$ is a limiting case of the
algebra ${\cal M}(R_{m|n}, R_{m|n})$ (\ref{REA}) at $q\to 1$ where
$R_{m|n}$ is the Drinfeld-Jimbo R-matrix obtained in quantization of
the classical supergroup $GL(m|n)$ (an explicit form of this
R-matrix can be found in \cite{DKS,I}). The so-called invariant
Cayley-Hamilton identity established in \cite{KT} is a limiting case
of the Cayley-Hamilton identity for the quantum matrix algebras of
the $GL(m|n)$ type. The latter is proved in the next section.

\section{Proof of the Cayley-Hamilton identities}
\label{sec:ch}

In this section we give a scheme of the proof and then formulate the
Cayley-Hamilton identity for the quantum matrix algebras
${\cal M}(R,F)$ defined by an R-matrix $R$ of the $GL(m|n)$ type.
Below these algebras will be referred to as {\em the matrix algebras
of the $GL(m|n)$ type}.

Fix a pair of integers $m\geq 1,\; n\geq 1$ and denote
$A:=(m+1)(n+1)$. We also introduce the special notations for the
following partitions
\ba
\lb{ch1}
\Lambda (r,s) &:=& \left((n+1)^r,n^{(m-r)},s\right)
, \quad
 r=0,\dots ,m,\; s=0,\dots ,n,
\\
\lb{ch2}
\Lambda^+ (r,s) &:=& \left(n+2,(n+1)^{(r-1)},n^{(m-r)},s\right),
\quad r=1,\dots ,m,\; s=0,\dots ,n,
\\
\lb{ch3}
\Lambda_+ (r,s) &:=& \left((n+1)^r,n^{(m-r)},s,1\right),
\quad  r=0,\dots ,m,\; s=1,\dots ,n.
\ea
Depicted below are the corresponding Young diagrams
$$
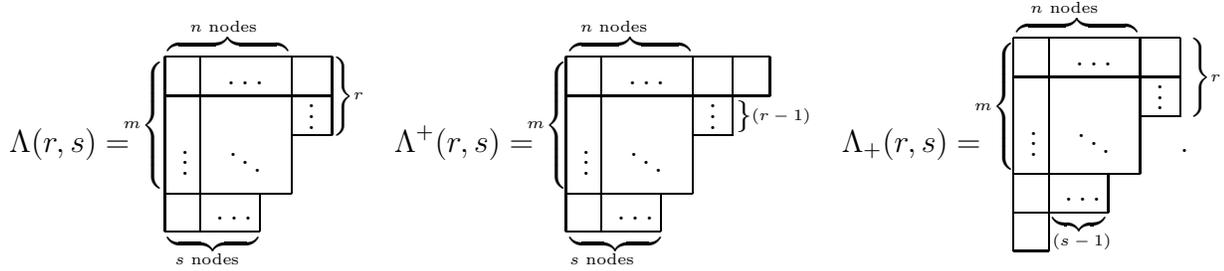
\lefteqn{
\hspace{24mm}\mbox{\raisebox{16mm}{\tiny $n$ nodes}}
\hspace{43mm}\mbox{\raisebox{16mm}{\tiny $n$ nodes}}
\hspace{51mm}\mbox{\raisebox{18.5mm}{\tiny $n$ nodes}}}
\lefteqn{
\hspace{21mm}\overbrace{\hspace{16mm}\rule{0pt}{13mm}}
\hspace{36.8mm}\overbrace{\hspace{16mm}\rule{0pt}{13mm}}
\hspace{42.7mm}\overbrace{\hspace{16mm}\rule{0pt}{15.5mm}}}
\lefteqn{
\hspace{22mm}\mbox{\raisebox{-15mm}{\tiny $s$ nodes}}
\hspace{44mm}\mbox{\raisebox{-15mm}{\tiny $s$ nodes}}
\hspace{55.7mm}\mbox{\raisebox{-12.5mm}{\tiny $(s-1)$ }}}
\lefteqn{
\hspace{21mm}\mbox{\raisebox{-11mm}{$\underbrace{\hspace{12mm}}$}}
\hspace{41.5mm}\mbox{\raisebox{-11mm}{$\underbrace{\hspace{12mm}}$}}
\hspace{52.2mm}\mbox{\raisebox{-8.4mm}{$\underbrace{
\hspace{7mm}}$}}}
\Lambda(r,s)=
\hspace{-2mm}\raisebox{3.2mm}{\tiny $m$}
\raisebox{2.5mm}{$\left\{\rule{0pt}{10mm} \right.$}
\hspace{-1pt}
\begin{tabular}{|c|cc|c|}\hline
&\multicolumn{2}{c|}{\hspace{-5mm}$\lefteqn\dots$}&\\ \hline
&&&\raisebox{-2pt}{$\smash{\vdots}$}\\ \cline{4-4}
\vdots$\!$&\multicolumn{2}{c|}{\hspace{-5mm}$\lefteqn\ddots$}&
\multicolumn{1}{c}{}\\[-12pt]
&\multicolumn{2}{c|}{}&\multicolumn{1}{c}{}\\ \cline{1-3}
&\multicolumn{1}{c|}{$\lefteqn\dots$ $\;\;$}&\multicolumn{2}{c}{}\\
\cline{1-2}
\end{tabular}
\hspace{-1pt}
\raisebox{6.3mm}{$\left. \rule{0pt}{6mm} \right\}$}
\raisebox{7mm}{\tiny $r$}
\quad \Lambda^+(r,s)=
\hspace{-2mm}\raisebox{3.2mm}{\tiny $m$}
\raisebox{2.5mm}{$\left\{\rule{0pt}{10mm} \right.$}
\hspace{-2pt}
\begin{tabular}{|c|cc|c|c|}\hline
&\multicolumn{2}{c|}{\hspace{-5mm}$\lefteqn\dots$}&&\\ \hline
&&&\raisebox{-2pt}{$\smash{\vdots}$}\hspace{-1pt}&
\multicolumn{1}{c}{\,}\\
\cline{4-4}
\vdots $\!$&\multicolumn{2}{c|}{\hspace{-5mm}$\lefteqn\ddots$}&
\multicolumn{2}{c}{}\\[-12pt]
&\multicolumn{2}{c|}{}&\multicolumn{2}{c}{}\\ \cline{1-3}
&\multicolumn{1}{c|}{$\lefteqn\dots$ $\;\;$}&\multicolumn{3}{c}{}\\
\cline{1-2}
\end{tabular}
\hspace{-14pt}
\raisebox{3.7mm}{$\left. \rule{0pt}{3mm} \right\}$}
\raisebox{4.2mm}{\tiny $(r-1)$}
\quad \Lambda_+(r,s)=
\hspace{-2mm}\raisebox{5.5mm}{\tiny $m$}
\raisebox{5mm}{$\left\{\rule{0pt}{10mm} \right.$}
\hspace{-2pt}
\begin{tabular}{|c|cc|c|}\hline
&\multicolumn{2}{c|}{\hspace{-5mm}$\lefteqn\dots$}&\\ \hline
&&&\raisebox{-2pt}{$\smash{\vdots}$}\\ \cline{4-4}
\vdots $\!$&\multicolumn{2}{c|}{\hspace{-5mm}$\lefteqn\ddots$}&
\multicolumn{1}{c}{}\\[-12pt]
&\multicolumn{2}{c|}{}&\multicolumn{1}{c}{}\\ \cline{1-3}
&\multicolumn{1}{c|}{$\lefteqn\dots$ $\;\;$}&\multicolumn{2}{c}{}\\
\cline{1-2}
&\multicolumn{3}{c}{}\\ \cline{1-1}
\end{tabular}.
\hspace{-1pt}
\raisebox{8.8mm}{$\left. \rule{0pt}{6mm} \right\}$}
\raisebox{9.3mm}{\tiny $r$}
$$

Let the symbol $E^{\Lambda(r,s)}_{row}$, $s\geq 1$, (respectively,
$E^{\Lambda(r,s)}_{col}$, $r\geq 1$, $E^{\Lambda^+(r,s)}$,
$E^{\Lambda_+(r,s)}$) denotes a diagonal matrix unit which
corresponds to a standard tableau of the form $\Lambda(r,s)$
(respectively, $\Lambda(r,s)$, $\Lambda^+(r,s)$, $\Lambda_+(r,s)$)
containing a standard Young tableau of the form $\Lambda(r,s-1)$
(respectively, $\Lambda(r-1,s)$, $\Lambda(r,s)$, $\Lambda(r,s)$).
In other words, in the diagonal matrix units
$E^{\Lambda(r,s)}_{row}$, $E^{\Lambda(r,s)}_{col}$,
$E^{\Lambda^+(r,s)}$ and $E^{\Lambda_+(r,s)}$ one fixes the position
of the largest number $S$ ($S$ is equal to $mn+r+s$, $mn+r+s$,
$mn+r+s+1$ and $mn+r+s+1$ respectively) as is shown in the picture
below
\be
\lb{ch4}
E^{\Lambda(r,s)}_{row} \!=\!
\begin{tabular}{|ccc|c|}\hline
\multicolumn{3}{|c|}{}&\vdots\\ \cline{4-4}
\raisebox{2pt}{$\smash{\lefteqn{\hspace{3mm}\ddots}}$}
&&&\multicolumn{1}{c}{}\\[-5pt]
\multicolumn{3}{|c|}{}&\multicolumn{1}{c}{}\\ \cline{1-3}
\multicolumn{1}{|c|}{\hspace{-5pt}$\lefteqn{\dots}\hspace{9pt}$}&
\multicolumn{1}{c|}{$\lefteqn{\hspace{-1pt}\scriptstyle S}
\hspace{2pt}$}& \multicolumn{2}{c}{}\\ \cline{1-2}
\end{tabular}
\;\;
E^{\Lambda(r,s)}_{col} \!=\!
\begin{tabular}{|cc|c|}\hline
\multicolumn{2}{|c|}{}&
\vspace{2mm}\smash{\raisebox{-2pt}{\vdots}}\\[-5pt] \cline{3-3}
\raisebox{-3pt}{$\smash{\lefteqn{\ddots}}$}&&
\hspace{-4pt}$\lefteqn{\scriptstyle S}$\\ \cline{3-3}
\multicolumn{2}{|c|}{}&\multicolumn{1}{c}{}\\ \cline{1-2}
\multicolumn{1}{|c|}{\hspace{-5pt}$\lefteqn{\dots}\hspace{9pt}$}&
\multicolumn{2}{c}{}\\ \cline{1-1}
\end{tabular}
\;\;
E^{\Lambda^+(r,s)}\! =\!
\begin{tabular}{|cc|c|c|}\hline
\multicolumn{2}{|c|}{}&&
\hspace{-2pt}$\hspace{1pt}\lefteqn{\scriptstyle S}\hspace{5pt}$\\
\cline{3-4}
\raisebox{-3pt}{$\smash{\lefteqn{\ddots}}$}&&
\raisebox{-2pt}{\smash{\vdots}}&
\multicolumn{1}{c}{}\\ \cline{3-3}
\multicolumn{2}{|c|}{}&\multicolumn{2}{c}{}\\ \cline{1-2}
\multicolumn{1}{|c|}{\hspace{-5pt}$\lefteqn{\dots}\hspace{9pt}$}&
\multicolumn{3}{c}{}\\ \cline{1-1}
\end{tabular}
\;\;
E^{\Lambda_+(r,s)} \!=\!
\begin{tabular}{|ccc|c|}\hline
\multicolumn{3}{|c|}{}&\vdots\\ \cline{4-4}
\raisebox{6pt}{$\smash{\lefteqn{\hspace{2mm}\ddots}}$}&&&
\multicolumn{1}{c}{}\\[-14pt]
\multicolumn{3}{|c|}{}&\multicolumn{1}{c}{}\\ \cline{1-3}
\multicolumn{2}{|c|}{\hspace{-5pt}
\raisebox{1pt}{$\lefteqn{\dots}\hspace{9pt}$}}&
\multicolumn{2}{c}{}\\ \cline{1-2}
\multicolumn{1}{|c|}{$\lefteqn{\hspace{-1pt}\scriptstyle S}
\hspace{2pt}$}
&\multicolumn{3}{c}{}\\ \cline{1-1}
\end{tabular}.
\ee
Emphasize, that the further calculations do not depend on the
positions of the other numbers in the nodes of the matrix units
(\ref{ch4}). At last, introduce the special notations for the
following elements of the space ${\rm Mat}_N({\cal M}(R,F))$
\ba
\lb{ch5}
P_{row}(r,s) &:=&
\Tr{2\dots A}\left(M_{\overline 2}\dots M_{\overline A}\
\rho_R(E^{\Lambda(r,s)}_{row})\, R_{t} R_{t-1}\dots R_1\right) ,
\\[2pt]
\lb{ch6}
P_{col}(r,s) &:=&
\Tr{2\dots A}\left(M_{\overline 2}\dots M_{\overline A}\
\rho_R(E^{\Lambda(r,s)}_{col})\, R_{t} R_{t-1}\dots R_1\right) ,
\\[2mm]
\lb{ch7}
P^+(r,s) &:=&
\Tr{2\dots A}\left(M_{\overline 2}\dots M_{\overline A}\
\rho_R(E^{\Lambda^+(r,s)})\, R_{t-1} R_{t-2}\dots R_1\right) ,
\\[2mm]
\lb{ch8}
P_+(r,s) &:=&
\Tr{2\dots A}\left(M_{\overline 2}\dots M_{\overline A}\
\rho_R(E^{\Lambda_+(r,s)})\, R_{t-1} R_{t-2}\dots R_1\right) ,
\ea
where $t:=(m-r)+(n-s)+1$.

Let us take the following linear combinations of the above matrices
\ba
\nonumber
\Phi_i &:=&
\sum_{k=\max\{0,i-n\}}^{\min\{i-1,m\}}
(-1)^{k}\,{(i-k)_q (m+n-i+k+2)_q\over (m+n-i+2)_q}\, P_{row}(k,i-k)
\\
&&\hspace{9mm}-\
\sum_{k=\max\{1,i-n\}}^{\min\{i,m\}}
(-1)^{k}\,{k_q (m+n-k+2)_q\over (m+n-i+2)_q}\, P_{col}(k,i-k)\ ,
\lb{ch9}
\ea
where the index $i$ runs from 1 to $(m+n)$. Our immediate task is
in proving the relations
\be
\lb{ch10}
\Phi_{i+1} - \Phi_i \ = \ \phi(M^{\overline{m+n-i}})
\sum_{k=\max\{0,i-n\}}^{\min\{i,m\}} (-1)^k\, q^{2k-i}\,
s_{\Lambda(k,i-k)}(M)\, ,\quad i=1,2,\dots ,m+n-1.
\ee
With that end in view, transform the expression of $\Phi_i$ to the
form
\ba
\nonumber
\Phi_i &=&{\textstyle
\sum\limits_{k=\max\{0,i-n\}}^{\min\{i-1,m\}}
 (-1)^{k}\,{(i-k)_q (m+n-i+k+2)_q\over (m+n-i+2)_q}
\left\{
{q^{(m+n+k-i+2)}\over (m+n+k-i+2)_q}P^+(k,i-k)\right.}\hspace{40mm}
\\[0mm]
\nonumber
&&\hspace{10mm} \textstyle{\left.
- {q^{-(i-k)}\over (i-k)_q}P_+(k,i-k)
+ q P_{row}(k,i-k+1)
+ {q^{(m+n-i+1)}\over (m+n-i+1)_q}P_{col}(k+1,i-k)\right\}}
\\[2mm]
\nonumber
&-&
\textstyle{
\sum\limits_{k=\max\{1,i-n\}}^{\min\{i,m\}}
(-1)^{k}\,{k_q (m+n-k+2)_q\over (m+n-i+2)_q}
\left\{
{q^k\over k_q}P^+(k,i-k)
-{q^{-(m+n-k+2)}\over (m+n-k+2)_q}P_+(k,i-k)\right.}
\\[0mm]
\lb{ch11}
&&\hspace{10mm}
\textstyle{
\left.
-{q^{-(m+n-i+1)}\over (m+n-i+1)_q}P_{row}(k,i-k+1)
- q^{-1}P_{col}(k+1,i-k)
\right\} .}
\ea
Here we extend the definition of elements (\ref{ch5})-(\ref{ch8})
to the boundary values of its indices, that is to the cases when
the matrix units used in construction of (\ref{ch5})-(\ref{ch8})
are not defined:
$$
P_{row}(k,n+1)=P_{col}(m+1,k)=P^+(k,0)=P_+(0,k)=0, \quad
k=0,1,\dots .
$$
When deriving (\ref{ch11}) we apply formulae (\ref{h-emb3}) to the
matrix units (\ref{ch4})
\be
\lb{ch12}
E^{\Lambda(r,s)}_{row/col}\, =\,
E^{\Lambda(r,s+1)}_{row} + E^{\Lambda(r+1,s)}_{col} +
E^{\Lambda^+(r,s)}+E^{\Lambda_+(r,s)}\, ,
\ee
and then use (\ref{me-con1a}) in the form
\be
\lb{ch13}
\rho_R (E^\lambda_\alpha) R_k \, =\,
\omega^{-1}(\ell_k)\,
\rho_R(E^\lambda_{\alpha,\pi_k(\alpha)})
\, -\, {q^{-\ell_k}\over (\ell_k)_q}\,
\rho_R(E^\lambda_\alpha)\, , \qquad
\lambda\vdash(k+1), \qquad
\ell_k\equiv \mtbl{\lambda}{\alpha}{k}\,
\ee
in order to transform the terms
$\rho_R(E^{\Lambda(r,s)}_{row/col})\, R_{t}$ in expressions for
$P_{row/col}(k,i-k)$. Using calculations analogous to (\ref{nul}),
one can show that the term with the non-diagonal matrix unit
entering the right hand side of (\ref{ch13}) does not contribute
to (\ref{ch11}).

Note that expressions (\ref{ch12}) for $E^{\Lambda(r,s)}_{row}$ and
$E^{\Lambda(r,s)}_{col}$ are {\em formally} the same, since our
notations (\ref{ch4}) take into account only the position of the
largest integer among those filling the Young tableau. One should
keep in mind that the matrix units $E^{\Lambda(r,s)}_{row}$ and
$E^{\Lambda(r,s)}_{col}$ in the right hand side of (\ref{ch12})
differ in the position of the number preceding to the largest one.
This difference manifests itself in applying relation (\ref{ch13}).
\medskip

At the next step we reduce the similar terms in (\ref{ch11}):
\ba
\nonumber
\Phi_i &=&{\textstyle
- \sum\limits_{k=\max\{0,i-n\}}^{\min\{i,m\}}
(-1)^{k}\, q^{2k-i} \left\{ P^+(k,i-k)  + P_+(k,i-k) \right\}
}
\\[0mm]
\nonumber
&+&
\textstyle{
\sum\limits_{k=\max\{0,i+1-n\}}^{\min\{i,m\}}
{(-1)^{k}\over (m+n-i+2)_q}
\Bigl\{q^{-(m+n-i+1)}\, {k_q (m+n-k+2)_q\over (m+n-i+1)_q}
\Bigr.}
\\[0mm]
\nonumber
&&\hspace{10mm}
\textstyle{
\Bigl.
+\ q\,  (i-k)_q (m+n-i+k+2)_q \Bigr\}
P_{row}(k,i-k+1)}
\\[2mm]
\nonumber
&+&
\textstyle{
\sum\limits_{k=\max\{0,i-n\}}^{\min\{i,m-1\}}
{(-1)^{k}\over (m+n-i+2)_q}
\Bigl\{q^{(m+n-i+1)}\, {(i-k)_q (m+n-i+k+2)_q\over (m+n-i+1)_q}
\Bigr.}
\\[0mm]
\lb{ch14}
&&\hspace{10mm}
\textstyle{
\Bigl.
+\ q^{-1} k_q (m+n-k+2)_q \Bigr\}
P_{col}(k+1,i-k) .}
\ea
In calculating the coefficients at $P^+(\cdot,\cdot)$ and
$P_+(\cdot,\cdot)$ the following relation is useful
\be
\lb{ch15}
q^x y_q + q^{-y} x_q = (x+y)_q, \quad x, y\in
{\Bbb Z}\hspace{-5pt}{\Bbb Z}.
\ee

Next, we transform the first summand in the right hand side of
(\ref{ch14}) with the help of the computation
\ba
\nonumber
&\lefteqn{\hspace{-10mm}
P^+(k,i-k)+P_+(k,i-k)+P_{row}(k,i-k+1)+P_{col}(k+1,i-k)
}&
\\[2mm]
\lb{ch16}
&=&
\Tr{2\dots A}\left(M_{\overline 2}\dots M_{\overline A}\
\rho_R(E^{\Lambda(r,s)}_{\alpha})\, R_{m+n-i}\dots R_1\right)
\\[2pt]
\nonumber
&=&
\Tr{2\dots (m+n-i+1)}\left(M_{\overline 2}\dots
M_{\overline{m+n-i+1}}\,
R_{m+n-i}\dots R_1\right)\, s_{\Lambda(k,i-k)}(M)
\\[2pt]
\lb{ch17}
&=&
\phi(M^{\overline{m+n-i}})\, s_{\Lambda(k,i-k)}(M)\, .
\ea
In the above computation we use relation (\ref{ch12}) first what
results in (\ref{ch16}). Expression (\ref{ch16}) turns out to be
independent on a particular choice of the index $\alpha$ of the
matrix unit. Then we split (\ref{ch16}) into two factors with the
use of (\ref{char1}).  At last, the first factor is transformed in
the same way as that used to derive relation (\ref{m-k2}) for the
matrix powers from definition (\ref{m-k}).

Substituting (\ref{ch17}) into (\ref{ch14}) one gets
\ba
\nonumber
\Phi_i
&=& -\ \phi(M^{\overline{m+n-i}})\,
\sum_{k=max\{0,i-n\}}^{\min\{i,m\}} (-1)^{k} q^{2k-i}\,
s_{\Lambda(k,i-k)}(M)
\\[2mm]
\nonumber
&&+\
\sum_{k=\max\{0,i+1-n\}}^{\min\{i,m\}}
(-1)^{k}\,{(i+1-k)_q (m+n-i+k+1)_q\over (m+n-i+1)_q}\,
P_{row}(k,i+1-k)
\\[2mm]
&&+\
\sum_{k=\max\{0,i-n\}}^{\min\{i,m-1\}}
(-1)^{k}\,{(k+1)_q (m+n-k+1)_q\over (m+n-i+1)_q}\,
P_{col}(k+1,i-k)\ .
\lb{ch18}
\ea
Here to simplify the coefficients at $P_{row}(\cdot ,\cdot)$ and
$P_{col}(\cdot ,\cdot)$ we have used the formula
\be
\lb{ch19}
\textstyle
q^{-x}\, {(x+y+1)_q z_q\over x_q} + q\, y_q (x+z+1)_q +
q^{z-y}\, (x+1)_q \ = \ {(x+1)_q (y+1)_q (x+z)_q\over x_q} ,
\ee
which can be easily verified with the help of (\ref{ch15}).

Finally, shifting the summation index $k\rightarrow k+1$ in the last
summand of (\ref{ch18}) we come to (\ref{ch10}).
\medskip

Now it remains only to construct analogs of (\ref{ch10}) for the
boundary values $i=0$ and $i=m+n$.

For $i=0$ we have:
\be
\lb{ch20}
\Phi_1 = P_{row}(1,0) + P_{col}(0,1) =
\phi(M^{\overline{m+n}})\, s_{\Lambda(0,0)}(M)\, ,
\ee
which coincides with (\ref{ch10}) if one sets $\Phi_0 = 0$.

In case $i=m+n$ the transformation of $\Phi_{m+n}$ is performed in
the same way as when deriving (\ref{ch10}). A peculiarity of the
case involved consists in two points. Firstly, when carrying out the
computations one finds an additional term associated with the Young
diagram $((n+1)^{m+1})$ and secondly, $\Phi_{m+n+1}=0$. So, making
the calculations we find
\ba
\nonumber
\Phi_{m+n}&  =&
(-1)^{m+1} q^{m-n}\, \mbox{\rm Id}_V\, s_{\Lambda(m,n)}(M)
\\[2mm]
\lb{ch21}
&+&
(-1)^{m}\, (m+1)_q (n+1)_q
\Tr{2\dots A}\Bigl\{M_{\overline 2}\dots M_{\overline A}\
\rho_R\left(E^{((n+1)^{m+1})}_{m+1}\right)\Bigr\} ,
\ea
which also has the same form as (\ref{ch10}) provided that
$\rho_R\left(E^{((n+1)^{m+1})}_{m+1}\right)=0$ or, by definition,
in the $GL(m|n)$ case.

Next, we add together the left and right hand sides of all relations
(\ref{ch10}), (\ref{ch20}) and then subtract from the sum the
corresponding parts of (\ref{ch21}). At last, on applying the map
$\phi^{-1}$ (see (\ref{fi}), (\ref{fi-inv})) to the result we come to
the following theorem.

\begin{teor}(The Cayley-Hamilton identity)
Let $M$ be the matrix of generators of a quantum matrix algebra
${\cal M}(R,F)$ of the $GL(m|n)$ type. Then the following matrix
identity holds
\be
\lb{ch22}
\sum_{i=0}^{n+m}  M^{\overline{ m+n-i}}\,
\sum_{k=\max\{0,i-n\}}^{\min\{i,m\}} (-1)^k\, q^{2k-i}\,
s_{\Lambda(k,i-k)}(M)
\,\equiv 0\, ,
\ee
\end{teor}

If one multiplies the both sides of the matrix equalities
(\ref{ch10}), (\ref{ch20}) and  (\ref{ch21}) by the matrix $M$ from
the left, the equalities turn into  relations among the elements of
the space ${\rm Pow}(R,F)$. Adding this relations together as in
deriving the Cayley-Hamilton theorem, we establish the connection
among the elements of two basis sets of the space ${\rm Pow}(R,F)$
which were described in proposition \ref{2bazisa}.
\begin{pred}
Let $M$ be the matrix of generators of a quantum matrix algebra
${\cal M}(R,F)$ of the Hecke type. Then the $\lambda$-powers of $M$
corresponding to the rectangular Young diagrams
$\lambda=\left((r+1)^{s+1}\right)$, $r,s=0,1,\dots$ are expressed in
terms of its $\overline{k}$-th powers as follows
\ba
\nonumber
\lefteqn{\hspace{-10mm}
(-1)^{s}\, (s+1)_q (r+1)_q
M^{\left(((r+1)^{s+1});\, s+1\right)}
}&&
\\[2mm]
\lb{ch23}
&=&
\sum_{i=0}^{s+r}  M^{\overline{ s+r+1-i}}\,
\sum_{k=\max\{0,i-r\}}^{\min\{i,s\}} (-1)^k\, q^{2k-i}\,
s_{\Lambda(k,i-k)}(M).
\ea
\end{pred}

\end{document}